\documentclass[a4paper,12pt]{article}

\usepackage{mathtools}
\usepackage{amsmath,amssymb,amsthm, bm}
\usepackage{enumerate}
\usepackage[noadjust]{cite}

\usepackage[spanish, english]{babel}

\usepackage[ansinew]{inputenc}
\usepackage{here}
\usepackage{hyperref}

\usepackage{xcolor}

\usepackage{graphicx}

\usepackage[normalem]{ulem}

\newcommand{\tp}{\operatorname{tp}}

\newenvironment{claimproof}[1][\proofname]
               {
                 \proof[#1]
                 
               }
               {
                 \endproof
               }

\newcommand{\FF}{\mathcal{F}}

\newcommand{\SSS}{\mathcal{S}}

\theoremstyle{definition}
\newtheorem{definition}{Definition}[section]

\theoremstyle{theorem}

\newtheorem{thmx}{Theorem}
 
\newtheorem{conx}[thmx]{Conjecture}
\newtheorem{corx}[thmx]{Corollary}
\newtheorem{lemma}[definition]{Lemma}
\newtheorem{proposition}[definition]{Proposition}

\newtheorem{claim}[definition]{Claim}

\newtheorem*{question*}{Question}

\newtheorem{questions}[definition]{Questions}

\theoremstyle{remark}
\newtheorem*{acknowledgment}{\textbf{Acknowledgments}}

\title{Definable $(\omega,2)$-theorem for families with VC-codensity less than $2$}
\author{Pablo And\'ujar Guerrero \\
\small \emph{Email address:} pa377@cantab.net}
\date{}

\begin{document}
\maketitle 

\noindent
{\small \emph{2020 Mathematics Subject Classification.} 03C45 (Primary), 52A35 (Secondary). \\
\emph{Key words.} NIP, VC-density, (p,q)-theorem.} 

\begin{abstract}
Let $\mathcal{S}$ be a family of nonempty sets with VC-codensity less than $2$. 
We prove that, if $\mathcal{S}$ has the $(\omega,2)$-property (for any infinitely many sets in $\mathcal{S}$, at least $2$ among them intersect), then $\mathcal{S}$ can be partitioned into finitely many subfamilies, each with the finite intersection property. If $\mathcal{S}$ is definable in some first-order structure, then these subfamilies can be chosen definable too. 

This is a strengthening of the case $q=2$ of the definable $(p,q)$-conjecture in model theory~\cite{simon15} and the Alon-Kleitman-Matou\v{s}ek $(p,q)$-theorem in combinatorics~\cite{matousek04}.

\end{abstract}

\section{Introduction}

Given a family of sets $\mathcal{S}$, a boolean atom is a maximal nonempty intersection of sets in the closure of $\mathcal{S}$ under complements. The dual shatter function $\pi^*_{\mathcal{S}}:\omega\rightarrow \omega$ of $\mathcal{S}$ sends each $n$ to the maximum number of boolean atoms of any subfamily of $\mathcal{S}$ of size $n$. 

For cardinals $p\geq q>1$, a family of sets $\mathcal{S}$ has the $(p,q)$-property if it does not contain the empty set and, for any $p$ sets in $\mathcal{S}$, there exists a subfamily among them of size $q$ with nonempty intersection.


 
Using ideas from Alon and Kleitman~\cite{alon-kleitman}, Matou\v{s}ek proved the following in~\cite[Theorem 4]{matousek04}. 

\begin{thmx}[Alon-Kleitman-Matou\v{s}ek $(p,q)$-theorem\footnote{While classically the Alon-Kleitman-Matou\v{s}ek $(p,q)$-theorem is stated for finite $\mathcal{F}$, a straightforward application of first-order logic compactness shows that this is equivalent to the infinite version presented here (see the proof of \cite[Proposition 2.5]{simon15}).}]\label{thm:matousek_2}
Let $q\geq 2$ be an integer and $\SSS$ be a family of sets whose dual shatter function satisfies $\pi_{\SSS}^*(n)\in o(n^q)$ (that is, \mbox{$\lim_{n\rightarrow \infty} \pi_{\SSS}^*(n)/n^q =0$}). For any integer $p\geq q$, there exists some $m<\omega$ such that, if $\FF$ is a subfamily of $\SSS$ with the $(p,q)$-property, then $\FF$ can be partitioned into at most $m$ subfamilies, each with the finite intersection property.  
\end{thmx}



For notational conventions and some model theoretic definitions in this paper we refer the reader to Section~\ref{section:definitions} and to~\cite{simon:guide-nip}.

Chernikov and Simon~\cite{cher_sim_15} used Theorem~\ref{thm:matousek_2} to study NIP theories. 
In~\cite[Problem 29]{cher_sim_15} they asked whether a definable version of it holds in this setting. This has evolved to be known as the definable $(p,q)$-conjecture \cite[Conjecture 2.15]{simon15}. Specifically, the conjecture (which was put forward before the connection with the $(p,q)$-theorem was established) states that any NIP formula which is non-dividing over a model $M$ belongs to a (finitely) consistent $M$-definable family. By means of first-order logic compactness, as well as Theorem~\ref{thm:matousek_2}, this can be restated as follows. 


\begin{conx}[Definable $(p,q)$-conjecture\footnote{In the literature the conjecture is commonly found with the stronger assumption that the whole structure is NIP~\cite[Conjecture 5.1]{simon15}. Kaplan~\cite[Corollary 4.9]{kaplan22} has recently presented a proof of this version of the conjecture.}]\label{conjecture:def_pq}
Let $q\geq 2$ be an integer, $M$ be an $L$-structure and $\varphi(x,y)$ be an $L(M)$-formula (which we identify with the family of sets $\{\varphi(M,a) : a\in M^{|y|}\}$) with dual shatter function $\pi_{\varphi}^*(n)\in o(n^q)$. If there exists an integer $p\geq q$ such that $\varphi(x,y)$ has the $(p,q)$-property, then there exists some $m<\omega$ and $L(M)$-formulas $\sigma_1(y),\ldots, \sigma_m(y)$ such that $\cup_i \sigma_i(M) = M^{|y|}$ and, for every $i\leq m$, the family $\{\varphi(x,a) : a\in \sigma_i(M)\}$ is consistent.    
\end{conx}

Conjecture~\ref{conjecture:def_pq}, which can be seen as a definable non-uniform version of Theorem~\ref{thm:matousek_2}, is known to hold in certain cases. Simon~\cite{simon14} proved it in dp-minimal theories for formulas $\varphi(x,y)$ with $|x|\leq 2$, and in any theory for formulas that extend to an invariant type of dp-rank $1$. 
In~\cite{simon15}, he proved it in NIP theories of small or medium directionality. Simon and Starchenko~\cite[Theorem 5]{simon_star_14} proved a stronger version of the conjecture for a class of dp-minimal theories that includes those that are linearly ordered, unpackable VC-minimal, or have definable Skolem functions. Recently, Boxall and Kestner~\cite{boxall_kestner_18} proved, using Theorem~\ref{thm:matousek_2} and the work on NIP forking of Chernikov and Kaplan~\cite{cher_kap_12}, Conjecture~\ref{conjecture:def_pq} in distal NIP theories. While this paper was under review, Itay Kaplan~\cite{kaplan22} presented a proof of a uniform version of Conjecture~\ref{conjecture:def_pq} for formulas in NIP theories.

In this paper we prove a strengthening of both Conjecture~\ref{conjecture:def_pq} and (the non-uniform version of) Theorem~\ref{thm:matousek_2} in the case where $q=2$. In particular, we show that Conjecture~\ref{conjecture:def_pq} holds when $q=2$, and that we may furthermore weaken the $(p,2)$-property to the $(\omega,2)$-property in the statements of Conjecture~\ref{conjecture:def_pq} and the case $\mathcal{S}=\mathcal{F}$ of Theorem~\ref{thm:matousek_2}. 

\begin{thmx}[Definable $(\omega,2)$-theorem]\label{them:intro_dfbl_pq}
Let $M$ be an $L$-structure and $\varphi(x,y)$ be an $L(M)$-formula with dual shatter function $\pi_{\varphi}^*(n)\in o(n^2)$ (e.g $VC$-codensity of $\varphi(x,y)$ is less than $2$). If $\varphi(x,y)$ has the $(\omega,2)$-property, then there exist some $m<\omega$ and $L(M)$-formulas $\sigma_1(y),\ldots, \sigma_m(y)$ such that $\cup_i \sigma_i(M) = M^{|y|}$ and, for every $i\leq m$, the family $\{ \varphi(x,a) : a\in \sigma_i(M)\}$ is consistent.  
\end{thmx}

Since any family of sets can be witnessed as a definable family in some structure, the following corollary is immediate. 

\begin{corx}[$(\omega,2)$-theorem]\label{cor:pq}
Let $\SSS$ be a family of sets with $\pi_{\SSS}^*(n)\in o(n^2)$. If $\SSS$ has the $(\omega,2)$-property, then it can be partitioned into finitely many subfamilies, each with the finite intersection property. 
\end{corx}

Our proof of Theorem~\ref{them:intro_dfbl_pq} is elementary in that it avoids the use of both the Alon-Kleitman-Matousek $(p,q)$-theorem (as well as its related fractional Helly theorem) and the work of Shelah, Simon and others on NIP theories.  

\begin{acknowledgment}
This research was supported by the Fields Institute for
Research in Mathematical Sciences, specifically by the 2021 Thematic Program on Trends in Pure and Applied Model Theory and the 2022 Thematic Program on Tame Geometry, Transseries and Applications to Analysis and Geometry. It was also supported by the EPSRC grant EP/V003291/1.

Pierre Simon provided some helpful comments regarding the main result of the paper. I moreover thank the anonymous referee for suggestions which greatly improved the presentation and for finding a mistake in an old version of Lemma~\ref{fact:types_p_i}. Lastly, I thank Pantelis Eleftheriou and Mervyn Tong for reading early versions of this paper and offering useful feedback.
\end{acknowledgment}

\section{Preliminaries}

\subsection{Notation} \label{section:definitions}

Throughout we fix two structures $M\preccurlyeq U$ in some language $L$, where $U$ realizes every type over $M$. For any $A\subseteq U$, let $L(A)$ denote the expansion of $L$ by formulas with parameters in $A$. 


Given a (partitioned) formula $\varphi(x,y)$, some $b\in U^{|y|}$ and $A\subseteq U^{|x|}$, let 
$\varphi(A,b)=\{a\in A : U\models \varphi(a,b)\}$. For $A\subseteq U$, we write $\varphi(A,b)$ instead of $\varphi(A^{|x|}, b)$. 
By ``definable set" we mean ``definable set in $M$ possibly with parameters", i.e. a set of the form $\varphi(M)$ for some $L(M)$-formula $\varphi(x)$.

We apply notions such as the $(p,q)$-property and dual shatter function to formulas $\varphi(x,y)$ by adopting the usual convention of identifying them with the family of sets $\{\varphi(M,a) : a\in M^{|y|}\}$. In the context of formulas, we refer to the finite intersection property as being (finitely) consistent, and to being pairwise disjoint as being pairwise inconsistent.


Given a formula $\varphi(x,y)$ and $A\subseteq U^{|y|}$, by a \emph{$\varphi$-type over $A$} we mean a maximal consistent collection $p(x)$ of formulas in $\{ \varphi(x,a), \neg\varphi(x,a) : a\in A\}$. 

Throughout, $n$,$m$, $i$, $j$, $k$ and $l$ are positive integers.  

\subsection{Preliminary results}

We present some preliminary lemmas on $\varphi$-types for formulas $\varphi(x,y)$ with $\pi^*_\varphi(n)\in o(n^2)$.

\begin{lemma}\label{lemma:vc_types}
Let $\varphi(x,y)$ be an $L(M)$-formula such that $\pi^*_\varphi(n)\in o(n^2)$. 
Suppose that there exists some $b\in U^{|y|}$ such that $\varphi(M,b)=\emptyset$. 
Then there exists $\theta(y)\in \tp(b/M)$ such that the elements of $\varphi(U,b)$ realize only finitely many $\varphi$-types over $\theta(M)$.  
\end{lemma}
\begin{proof}
Let $\varphi(x,y)$ and $b\in U^{|y|}$ be as in the lemma. We assume that, for any $\theta(y)\in \tp(b/M)$, the elements of $\varphi(U,b)$ realize infinitely many $\varphi$-types over $\theta(M)$. We prove the lemma by showing that, for every $n$, 
\begin{equation}\label{eqn:codensity}
\pi^*_\varphi(n) \geq \sum_{i=1}^{n} i = \frac{n^2+n}{2}. 
\end{equation}
In particular, it follows that $\pi^*_\varphi(n)\notin o(n^2)$. 

We construct a sequence $(a_n : 1 \leq n<\omega)$ in $M^{|y|}$ and a set 
$\{ c_{i,j}:\, 1 \leq i<\omega,\, 1 \leq j \leq i\}$
in $M^{|x|}$ with the following property. 
For every $n$ and distinct pairs $(i,j)$, $(i',j')$, with $i, i'\leq n$, $j\leq i$ and $j'\leq i'$, it holds that
\begin{equation}\label{eqn:c_ij}
\varphi(c_{i,j}, \{a_1,\ldots, a_n\})\neq \varphi(c_{i',j'}, \{a_1,\ldots, a_n\}).
\end{equation}
That is, for every $n$, the set $\{ c_{i,j}: 1 \leq i\leq n, \, 1 \leq j\leq i\}$ witnesses that 
\[
|\{\varphi(c, \{a_1,\ldots, a_n\}) : c\in M^{|x|}\}|\geq \sum_{i=1}^{n} i,
\]
which in turn shows that the elements $\{a_1,\ldots, a_n\}$ witness Equation~\eqref{eqn:codensity}. 
Specifically, the set $\{c_{i,j} : 1 \leq i<\omega, 1 \leq j\leq i\}$ will have the following two properties:
\begin{enumerate}[(i)]
\item \label{itm:c(i,j)1}$\neg\varphi(c_{i',j'}, a_i) \text{ and }  \varphi(c_{i,j}, a_i) \text{ holds for all } i'<i, \, j'\leq i',\, j\leq i,$
\item \label{itm:c(i,j)2} $\varphi(c_{i,j}, \{a_1,\ldots, a_{i-1}\})\neq \varphi(c_{i,j'}, \{a_1,\ldots, a_{i-1}\}) \text{ for all } i\geq 2,\, j<j'\leq i.$
\end{enumerate}
It is easy to see that condition~\eqref{eqn:c_ij} follows from~\eqref{itm:c(i,j)1} and~\eqref{itm:c(i,j)2}.  

For every $n$ and $a_1,\ldots, a_n$ in $M^{|y|}$, let $s(a_1,\ldots, a_n)$ denote the number of boolean atoms $C$ of $\{\varphi(U,a_1),\ldots, \varphi(U,a_n)\}$ satisfying that \mbox{$\varphi(C,b)\neq \emptyset$}. We construct our sequence in such a way that $s(a_1,\ldots, a_n)\geq n+1$ for every $n$.

We proceed to build sets $\{ a_i : 1 \leq i\leq n\}$ and $\{c_{i,j} : 1 \leq i\leq n,\, 1 \leq j\leq i\}$ by induction on $n$.

\smallskip\noindent\textbf{Case $n=1$.} 

Since, by assumption, the elements of $\varphi(U,b)$ realize infinitely many $\varphi$-types over $M$, there must be some $a\in M^{|y|}$ such that
\[
\varphi(U,b)\cap \varphi(U,a)\neq \emptyset \text{ and } \varphi(U,b)\setminus \varphi(U,a)\neq \emptyset.
\]
Let $a_1$ be any such $a$. Let $c_{1,1}$ be any element in $\varphi(M,a_1)$. Observe that $s(a_1)=2$. 

\smallskip\noindent\textbf{Induction $n>1$.} 

Suppose we have a sequence $(a_1,\ldots, a_{n-1})$ in $M^{|y|}$ as desired. Since $s(a_1,\ldots, a_{n-1})\geq n$, there are $n$ distinct boolean atoms $C_1,\ldots, C_n$ of the family $\{\varphi(U,a_1),\ldots, \varphi(U,a_{n-1})\}$ containing each elements from $\varphi(U,b)$. Let 
\[
\theta(M)=\{ a\in M^{|y|} : \neg\varphi(c_{i,j}, a),\, \varphi(C_k,a)\neq\emptyset \text{ for } j\leq i < n, \, k\leq n\}.
\] 
Since $\varphi(M,b)=\emptyset$, note that $b\in \theta(U)$. Consequently, by assumption, the elements of $\varphi(U,b)$ realize infinitely many $\varphi$-types over $\theta(M)$. In particular, there must exist some boolean atom $C$ of $\{\varphi(U,a_1),\ldots, \varphi(U,a_{n-1})\}$ satisfying that the elements of $\varphi(C,b)$ realize more than one $\varphi$-type over $\theta(M)$. Let $a_n \in \theta(M)$ witness this, i.e. $\varphi(C,b)\cap \varphi(U,a_n)\neq \emptyset$ and $\varphi(C,b) \setminus\varphi(U,a_n)\neq \emptyset$. It then follows that $s(a_1,\ldots, a_n)\geq n+1$. 

Finally, by definition of $\theta(M)$, we have that $\varphi(C_j, a_n)\neq \emptyset$ for every $j\leq n$. For any $j \leq n$, let $c_{n,j}$ be an element in $\varphi(C_j, a_n) \cap M^{|x|}$. Then clearly $\{c_{i,j}: 1 \leq i\leq n,\, 1\leq j\leq i\}$ satisfies condition~\eqref{itm:c(i,j)2}.  By definition of $\theta(M)$, note that it also satisfies condition~\eqref{itm:c(i,j)1}. 
\end{proof}

\begin{lemma}\label{claim:thm_0}
Let $\varphi(x,y)$ be an $L(M)$-formula such that $\pi^*_\varphi(n)\in o(n^2)$. 
Suppose that there exists some $b\in U^{|y|}$ such that, for any $\sigma(y) \in \tp(b/M)$, the family $\{\varphi(x,a) : a\in \sigma(M)\}$ fails to be consistent. Then there exists $\theta(y)\in \tp(b/M)$ such that the elements of $\varphi(U,b)$ realize only finitely many
$\varphi$-types over $\theta(M)$ and, moreover,
for any such type $p(x)$ exactly one of the following two conditions holds.
\begin{enumerate}[(a)]
\item \label{eqn:theta_1} $\{ a\in \theta(M) : \varphi(x,a)\in p(x)\}=\emptyset$. 
\item \label{eqn:theta_2}
For every $\theta'(y)\in \tp(b/M)$, the set $\{ a \in \theta'(M)  : \varphi(x,a)\in p(x) \}$ is not definable (in $M$). 
\end{enumerate}
\end{lemma}
\begin{proof}
Note that, by definition of $b$, for any $c\in M^{|x|}$ we have $\varphi(c,y)\notin \tp(b/M)$. So $\varphi(M,b)=\emptyset$. We apply Lemma~\ref{lemma:vc_types}. Hence let $\theta_0(y)\in \tp(b/M)$ be such that the elements of $\varphi(U,b)$ realize only finitely many $\varphi$-types over $\theta_0(M)$. Since otherwise the lemma is trivial we may assume that $\varphi(U,b)\neq \emptyset$. We denote these types by $p_1(x),\ldots,p_m(x)$.

Let $F\subseteq \{1,\ldots, m\}$ be the set of $i$ satisfying that
there exists a formula $\theta_i(y)\in \tp(b/M)$ such that the set $\sigma_i(M)=\{ a\in \theta_i(M) : \varphi(x,a) \in p_i(x)\}$ is definable.
Observe that, for any $i\in F$, since $\{\varphi(x,a) : a\in \sigma_i(M)\}$ is consistent, by definition of $b$ it holds that $b\notin \sigma_i(M)$. Finally let $\theta(y)$ be given by
\[
\theta_0(y) \wedge \bigwedge_{i\in F} (\theta_i(y) \wedge \neg\sigma_i(y)). 
\]
Since $\theta(M)\subseteq \theta_0(M)$, the $\varphi$-types over $\theta(M)$ realized in $\varphi(U,b)$ are exactly the restrictions $p_i(x)|_{\theta(M)}$ of the types $p_i(x)$ to $\theta(M)$, for $i\leq m$. We have ensured that, for any $i\in F$, the type $p_i(x)|_{\theta(M)}$ is the (necessarily unique) type described by condition~\eqref{eqn:theta_1}. On the other hand, by definition of $F$, for any $j \in \{1,\ldots,m\}\setminus F$ the type $p_j(x)|_{\theta(M)}$ satisfies condition~\eqref{eqn:theta_2}.
\end{proof}

\begin{lemma}\label{fact:types_p_i}
 Let $\varphi(x,y)$, $b\in U^{|y|}$, $\theta(y)\in \tp(b/M)$ and $p(x)$ be such that they satisfy condition~\eqref{eqn:theta_2} in Lemma~\ref{claim:thm_0}. Then, for any $L(M)$-formula $\lambda(x)$ satisfying that $\varphi(U,b)\subseteq \lambda(U)$, there exists some $a\in \theta(M)$ such that 
\[
 \varphi(U,a)\subseteq \lambda(U) \text{ and } \varphi(x,a)\in p(x).
\] 
\end{lemma}
\begin{proof}
Let $\theta'(M)$ be the set of $a\in \theta(M)$ with $\varphi(U,a)\subseteq \lambda(U)$. Observe that $\theta'(y)\in \tp(b/M)$. Then, by condition~\eqref{eqn:theta_2} in Lemma~\ref{claim:thm_0}, the set $\{ a \in \theta'(M)  : \varphi(x,a)\in p(x) \}$ is nonempty. Let $a$ be any element in the set.
\end{proof}

\section{Proof of the main result}


We prove Theorem~\ref{them:intro_dfbl_pq} through the next proposition. 

\begin{proposition}\label{prop:thm_pq_1}
Let $\varphi(x,y)$ be an $L(M)$-formula with $\pi^*_\varphi(n)\in o(n^2)$ and suppose that there exists $b\in U^{|y|}$ such that, for any $\sigma(y) \in \tp(b/M)$, the family $\{\varphi(x,a) : a\in \sigma(M)\}$ fails to be consistent. 
Let $\chi(x)$ be an $L(M)$-formula such that $\varphi(U,b)\subseteq \chi(U)$. Then there exists some $a\in M^{|y|}$ such that 
\begin{equation*} 
\varphi(U,a)\subseteq \chi(U)
\end{equation*}
and moreover 
\begin{equation*} 
\varphi(U,a)\cap \varphi(U,b)=\emptyset. 
\end{equation*}
\end{proposition}
\begin{proof}
By Lemma~\ref{claim:thm_0}, there exists some $\theta(y)\in \tp(b/M)$ such that the elements of $\varphi(U,b)$ realize only finitely many $\varphi$-types over $\theta(M)$, and furthermore for any such type condition~\eqref{eqn:theta_1} or condition~\eqref{eqn:theta_2} in the lemma holds. By passing from $\theta(M)$ to $\theta(M) \cap \{ a \in M^{|y|} : \varphi(U,a)\subseteq \chi(U)\}$ if necessary, we may also assume that every $a\in\theta(M)$ satisfies that $\varphi(U,a)\subseteq \chi(U)$. In particular, to prove Proposition~\ref{prop:thm_pq_1} it suffices to find some $a\in \theta(M)$ such that $\varphi(U,a)\cap \varphi(U,b)=\emptyset$. Since otherwise the result is trivial we may assume that $\varphi(U,b)\neq \emptyset$.

Let $p_1(x), \ldots, p_l(x)$ denote the distinct $\varphi$-types over $\theta(M)$ realized by elements of $\varphi(U,b)$. We prove Proposition~\ref{prop:thm_pq_1} by finding some $a\in \theta(M)$ such that $\varphi(x,a)\notin p_i(x)$ for every $i\leq l$. If $l=1$ and $p_1(x)$ is the (unique) type described by condition~\eqref{eqn:theta_1} in Lemma~\ref{claim:thm_0}, then clearly it suffices to take any $a\in \theta(M)$ and we are done. We assume this is not the case.

Let the numbering of the types $p_i(x)$ be such that, for some fixed $k \in \{l-1, l\}$, the types $p_i(x)$ for $1\leq i \leq k$ satisfy condition~\eqref{eqn:theta_2} and the possibly remaining type $p_{i}(x)$ for $k<i\leq l$ satisfies condition~\eqref{eqn:theta_1} in Lemma~\ref{claim:thm_0}. 
Hence, either $k=l$ or otherwise $1\leq k=l-1$ and the type $p_l(x)$ satisfies that $\varphi(x,a)\notin p_l(x)$ for every $a\in \theta(M)$. In either case it suffices to find some $a\in\theta(M)$ with $\varphi(x,a)\notin p_i(x)$ for every $1\leq i \leq k$.

Now let us fix, for every $1 \leq i \leq k$, an $L(M)$-formula $\chi_i(x)$ satisfying the following conditions:
\begin{itemize}
\item $p_i(x)\models \chi_i(x)$ for every $i < k$,
\item $p_j(x) \models \chi_{k}(x)$ for all $k \leq j \leq l$,
\item $\chi_i(U)\cap \chi_j(U)=\emptyset$ for every $i<j\leq k$.
\end{itemize} 

We define, for any $1\leq m \leq k$ and elements $a_1,\ldots, a_{m-1} \in M^{|y|}$, a set $\psi_{m}(M, a_1,\ldots, a_{m-1}) \subseteq \theta(M)$ as follows.

For $m=k$, let $\psi_k(M, a_1,\ldots, a_{k-1})$ denote the set of all $a\in \theta(M)$ such that 
\[
\varphi(U,a)\subseteq \bigcup_{i=1}^{k-1} (\varphi(U, a_i)\cap \chi_i(U))\cup \chi_k(U).
\] 
For $m < k$, let $\psi_{m}(M, a_1,\ldots, a_{m-1})$ denote the set of all $a\in \theta(M)$ such that 
\[
\varphi(U,a)\subseteq \bigcup_{i=1}^{m-1} (\varphi(U,a_i)\cap \chi_i(U)) \cup \bigcup_{i=m}^k \chi_i(U)
\]
and moreover there exists two elements
$a', a''\in \psi_{m+1}(M, a_1,\ldots, a_{m-1}, a)$,
with
\[
\varphi(U,a')\cap \varphi(U,a'')\cap \chi_{m+1}(U)=\emptyset.
\]
\begin{claim}\label{claim:psi-dfbl}
For any $m \leq k$, the sets 
$
\psi_{m}(M, a_1,\ldots, a_{m-1})
$
are definable uniformly (in $M$) over the parameters $a_i \in M^{|y|}$, $i<m$. 
\end{claim}
\begin{claimproof}
For any given $m\leq k$, let $(\text{A}_m)$ be the statement that the sets $\psi_{m}(M, a_1,\ldots, a_{m-1})$
are definable uniformly over the parameters $a_i \in M^{|y|}$, $i<m$. Statement $(\text{A}_k)$ clearly holds by definition. 
Then, for any $m<k$, $(\text{A}_m)$ follows easily from $(\text{A}_{m+1})$ and the definition of sets $\psi_m(M, a_1,\ldots, a_{m-1})$.
\end{claimproof}

We now prove two claims regarding the set $\psi_{1}(M)$ that will yield Proposition~\ref{prop:thm_pq_1}, by showing the existence of some $a\in \theta(M)$ with $\varphi(x,a)\notin p_i(x)$ for every $i\leq k$.

\begin{claim}\label{claim:thm_pq_2}
There exist $a, a'\in \psi_1(M)$ 
such that
\[
\varphi(U,a)\cap \varphi(U,a')\cap \chi_{1}(U)=\emptyset.
\] 
\end{claim}
\begin{claimproof}
For any $m \leq k$ consider the following two statements $(\text{I}_m)$ and $(\text{II}_m)$:
\begin{description}
\item[$(\text{I}_m)$] Let $a_i \in M^{|y|}$ be such that $\varphi(x,a_i) \in p_i(x)$, for $i < m$, and let $a\in \theta(M)$. Suppose that
\[
\varphi(U,a)\subseteq \bigcup_{i=1}^{m-1} (\varphi(U, a_i)\cap \chi_i(U)) \cup \bigcup_{i=m}^k \chi_i(U)
\]
and
\[
\varphi(x,a)\in p_m(x).
\] 
Then 
\[
a\in \psi_m(M, a_1,\ldots, a_{m-1}).
\]
\item[$(\text{II}_m)$] Let $a_i \in M^{|y|}$ be such that $\varphi(x,a_i) \in p_i(x)$, for $i < m$. Then there exist 
\[
a, a'\in \psi_m(M, a_1,\ldots, a_{m-1})
\] 
such that
\[
\varphi(U,a)\cap \varphi(U,a')\cap \chi_{m}(U)=\emptyset.
\] 
\end{description}

We prove $(\text{I}_m)$ and $(\text{II}_m)$ for every $m\leq k$ using a reverse induction on $m$. Claim~\ref{claim:thm_pq_2} is then given by $(\text{II}_1)$.

Trivially $(\text{I}_k)$ holds by definition of $\psi_k(M, a_1,\ldots, a_{k-1})$, even without the condition $\varphi(x,a)\in p_k(x)$.
We prove the remaining statements as follows. For $m\leq k$, we derive $(\text{II}_m)$ from $(\text{I}_m)$ using Claim~\ref{claim:psi-dfbl}. For $m<k$, we derive $(\text{I}_m)$ from $(\text{II}_{m+1})$.

\smallskip\noindent\textbf{Proof of $(\text{I}_m)\Rightarrow (\text{II}_{m})$ for $m\leq k$.}

Let $\varphi(x,a_i)\in p_i(x)$ for $i < m$. 
Let $\theta'(M)$ be the set of all $a\in \theta(M)$ such that 
\[
\varphi(U,a)\subseteq \bigcup_{i=1}^{m-1} (\varphi(U,a_i)\cap \chi_i(U)) \cup \bigcup_{i=m}^k \chi_i(U).
\]
Note that $\theta'(y)\in \tp(b/M)$. By definition of $p_m(x)$ (see condition~\eqref{eqn:theta_2} in Lemma~\ref{claim:thm_0}), the set $A$ of all $a\in \theta'(M)$ with
$
\varphi(x,a)\in p_m(x)
$ 
is not definable (in $M$). By $(\text{I}_m)$ note that 
\[
A \subseteq \psi_m(M, a_1,\ldots, a_{m-1}).
\]
By Claim~\ref{claim:psi-dfbl}, the set $\psi_m(M, a_1,\ldots, a_{m-1})$ is definable. Since the subset $A$ is not definable, there must exist some $a\in \psi_m(M, a_1,\ldots, a_{m-1})$ that is not in $A$, in particular 
\[
\varphi(x,a)\notin p_m(x). 
\]

Now, by Lemma~\ref{fact:types_p_i}, there exists some $a'\in \theta(M)$ with 
\[
\varphi(U,a')\subseteq \bigcup_{i=1}^{m-1} (\varphi(U,a_i) \cap \chi_i(U)) \cup (\chi_m(U)\setminus \varphi(U,a)) \cup \bigcup_{i=m+1}^k \chi_i(U)
\]
such that 
\[
\varphi(x,a')\in p_m(x).
\]
(In the case $m=k=l-1$ Lemma~\ref{fact:types_p_i} can still be applied because $\varphi(x,a)\notin p_l(x)$ by definition of the type $p_l(x)$.) 
Once again by $(\text{I}_m)$ it follows that
\[
a'\in \psi_m(M, a_1,\ldots, a_{m-1}).
\]

Finally, by construction note that 
\[
\varphi(U,a)\cap \varphi(U,a') \cap \chi_m(U) = \emptyset. 
\]

\smallskip\noindent\textbf{Proof of $(\text{II}_{m+1})\Rightarrow (\text{I}_{m})$ for $m<k$.}

Let $\varphi(x,a_i)\in p_i(x)$ for $i < m$, and $a\in \theta(M)$ be as described in $(\text{I}_m)$. In particular we have that $\varphi(x,a)\in p_m(x)$. 

By $(\text{II}_{m+1})$, there exist 
$
a', a''\in \psi_{m+1}(M, a_1,\ldots, a_{m-1}, a)
$
such that 
\[
\varphi(U,a')\cap \varphi(U,a'')\cap \chi_{m+1}(U)=\emptyset.
\] 
But then by definition this means that  
$a\in \psi_m(M, a_1,\ldots, a_{m-1})$.
\end{claimproof}

\begin{claim}\label{claim:thm_pq_2.1}
Suppose that there exists some $a'\in \psi_1(M)$ with 
\[
\varphi(x,a')\notin p_1(x).
\]
Then there exists some $a\in \theta(M)$ satisfying that
\[
\varphi(x,a)\notin p_i(x) \text{ for every } 1 \leq i \leq k. 
\]
\end{claim}
\begin{claimproof}
For any $m \leq k$ consider the following statement $(\text{B}_m)$:
\begin{description}
\item[$(\text{B}_m)$] Let $a_i \in M^{|y|}$, $i<m$, be such that there exist
$
a'\in \psi_m(M, a_1,\ldots, a_{m-1}),
$
with 
\[
\varphi(x,a')\notin p_m(x).
\]
Then there exists some $a\in \theta(M)$ with 
\[
\varphi(U,a)\subseteq \bigcup_{i=1}^{m-1} \left(\varphi(U,a_i)\cap \chi_i(U)\right) \cup \bigcup_{i=m}^k \chi_i(U)
\]
satisfying that
\[
\varphi(x,a)\notin p_j(x) \text{ for every } m \leq j \leq k. 
\]
\end{description}
We prove $(\text{B}_{m})$ for every $m \leq k$ by reverse induction on $m$. Claim~\ref{claim:thm_pq_2.1} then immediately follows from $(\text{B}_1)$.
Let $a_i$, for $i<m$, and $a'$ be as in $(\text{B}_{m})$. 

For the base case $m=k$, it clearly suffices to take $a=a'$.
We assume that $m<k$ and show that $(\text{B}_{m+1}) \Rightarrow (\text{B}_{m})$. 

By definition of $\psi_m(M, a_1,\ldots, a_{m-1})$, there exist $a'', a''' \in \psi_{m+1}(M, a_1,\ldots, a_{m-1}, a')$ with
\[
\varphi(U,a'') \cap \varphi(U,a''') \cap \chi_{m+1}(U)=\emptyset.
\]
Without loss of generality we may assume that $\varphi(x,a'')\notin p_{m+1}(x)$.
By $(\text{B}_{m+1})$, we derive that there exists some $a\in \theta(M)$ such that 
\begin{equation}\label{eqn:claim2}
\varphi(U,a)\subseteq \bigcup_{i=1}^{m-1} \left(\varphi(U,a_i)\cap \chi_i(U)\right) \cup (\varphi(U,a')\cap \chi_m(U)) \cup \bigcup_{i=m+1}^k \chi_i(U)
\end{equation}
and 
\[
\varphi(x,a)\notin p_j(x) \text{ for every } m < j \leq k. 
\]
However, since $\varphi(x,a')\notin p_m(x)$, then by~\eqref{eqn:claim2} it must also be that $\varphi(x,a)\notin p_m(x)$.
\end{claimproof}

We now complete the proof of the proposition. By Claim~\ref{claim:thm_pq_2}, let $a', a''\in \psi_1(M)$ be two elements such that $\varphi(U,a') \cap \varphi(U,a'') \cap \chi_1(U)=\emptyset$. Without loss of generality we may assume that $a'$ is such that $\varphi(x,a')\notin p_1(x)$. By Claim~\ref{claim:thm_pq_2.1} we conclude that there exists some $a\in \theta(M)$ satisfying that $\varphi(x,a)\notin p_i$ for every $i\leq k$, as desired.
\end{proof}


\begin{proof}[Proof of Theorem~\ref{them:intro_dfbl_pq}]
Let $\varphi(x,y)$ be an $L(M)$-formula with $\pi^*_\varphi(n)\in o(n^2)$. We assume that $\varphi(x,y)$ does not partition into finitely many consistent families and derive that it does not have the $(\omega,2)$-property, i.e. we build a sequence $(a_n : 1\leq n<\omega)$ in $M^{|y|}$ such that the family $\{ \varphi(x,a_n) : 1 \leq n<\omega\}$ is pairwise inconsistent. 

Hence we assume that $\varphi(x,y)$ satisfies that, for any finite collection of $L(M)$-formulas $\{ \sigma_i(y) : 1 \leq i\leq m\}$, if the family $\{\varphi(x,a) : a\in \sigma_i(M)\}$ is consistent for every $i\leq m$, then there exists some $a\in M^{|y|}$ such that $a\notin \cup_i \sigma_i(M)$. By model theoretic compactness we may fix some $b\in U^{|y|}$ satisfying that, for any formula $\sigma(y)\in\tp(b/M)$, the family $\{\varphi(x,a) : a\in \sigma(M)\}$ fails to be consistent.
We build our sequence $(a_n : 1 \leq n<\omega)$ using Proposition~\ref{prop:thm_pq_1}. In particular it will satisfy that, for every $i<\omega$, it holds that
\begin{equation}\label{eqn:pf-thmC}
\varphi(U,a_i)\cap \varphi(U,b)=\emptyset
\end{equation}
We proceed inductively on $n$. 

By Proposition~\ref{prop:thm_pq_1} (with $\chi(x):=``x=x"$), let $a_1\in M^{|y|}$ be any element satisfying~\eqref{eqn:pf-thmC}. Then, for the inductive step, let $(a_1,\ldots, a_{n-1})$ be elements each satisfying~\eqref{eqn:pf-thmC} and such that the formulas $\varphi(x,a_i)$, for $i < n$, are pairwise inconsistent. Let $\chi(x)$ denote the formula
\[
\bigwedge_{i=1}^{n-1} \neg\varphi(x,a_i). 
\]
Note that $\varphi(U,b)\subseteq \chi(U)$. Now, applying Proposition~\ref{prop:thm_pq_1}, let $a_n\in M^{|y|}$ be an element satisfying~\eqref{eqn:pf-thmC} and $\varphi(U,a_n)\subseteq \chi(U)$. The family $\{ \varphi(x,a_i) : 1 \leq i\leq n\}$ is pairwise inconsistent as desired.
\end{proof}

We end the paper with some questions. We note that, while this paper was under review, Kaplan~\cite{kaplan22} presented a positive answer to Question (2) for formulas in NIP theories. 

\begin{questions} 
\quad

\begin{enumerate}[(1)]
\item \emph{Definable $(\omega,q)$-conjecture:}  Let $\varphi(x,y)$ be a formula and $q \geq 2$ an integer such that $\pi^*_\varphi(n)\in o(n^q)$. If $\varphi(x,y)$ has the $(\omega,q)$-property, does it partition into finitely many consistent definable subfamilies?
\item \emph{Uniform definable $(p,2)$-conjecture 1:} Let $\varphi(x,y)$ and $\psi(y,z)$ be formulas where $\pi^*_\varphi(n)\in o(n^2)$. Given any integer $p\geq 2$, is there an $m$ such that any family of the form $\{ \varphi(x,a) : M\models \psi(a,b) \}$, for $b\in M^{|z|}$, with the $(p,2)$-property partitions into at most $m$ consistent definable subfamilies? 
\item \emph{Uniform definable $(p,2)$-conjecture 2:} Let $\varphi(x,y)$ be a formula with $\pi^*_\varphi(n)\in o(n^2)$. Given any integer $p\geq 2$, is there an $m$ such that any definable subfamily of $\varphi(x,y)$ with the $(p,2)$-property partitions into at most $m$ consistent definable subfamilies? 
\end{enumerate}
\end{questions}



\bibliography{mybib_pq_base_case}
\bibliographystyle{alpha}

\end{document}